\documentclass{article}
\usepackage{amsfonts,amsmath}
\mathchardef\mhyphen="2D 
\usepackage[a4paper]{geometry}
\usepackage{microtype}
\usepackage{enumerate} 
\usepackage{relsize}

\def\zo/{$0\mkern2mu\mhyphen1$}

\def\nn/{$n \times n$}

\title{A PROOF of Weak Graph Positivity, for a Large Range of the Parameters}
\date{\today} 
\author{Paul Federbush\\
Department of Mathematics\\
University of Michigan\\
Ann Arbor, MI, 48109-1043}

\usepackage{amsthm} 
\newtheorem{thm}{Theorem}[section]
\newtheorem{conj*}{Conjecture} 
\newtheorem{lemma}{Lemma}
\numberwithin{lemma}{section}
\numberwithin{conj*}{section}
\numberwithin{equation}{section} 
\DeclareMathOperator{\Prob}{Prob}
\DeclareMathOperator{\EE}{E}
\DeclareMathOperator{\OO}{\mathcal{O}}
\DeclareMathOperator{\PP}{P}
\newcommand{\sslash}{\mathbin{/\mkern-6mu/}} 
\newcommand{\stirlingi}{\genfrac{[}{]}{0pt}{}}

\newcommand{\stirlingii}{\genfrac{\{}{\}}{0pt}{0}}
\newcommand{\symm}{\mathrm{Sym}}

\begin{document}

\maketitle
\begin{abstract}
   One deals with $r$-regular bipartite graphs with $2n$ vertices.
  In a previous paper Butera, Pernici, and the author have introduced a quantity
  $d(i)$, a function of the number of $i$-matchings, and conjectured that as $n$
  goes to infinity the fraction of graphs that satisfy $\Delta^k d(i) \ge 0$,
  for all $k$ and $i$, approaches $1$.
  Here $\Delta$ is the finite difference operator. This conjecture we called the ``graph
  positivity conjecture''.
  ``Weak graph positivity'' is the conjecture that for each $i$ and $k$ the probability that $\Delta^kd(i)\geq0$ goes to $1$ as $n$ goes to infinity. Here we prove this for the range of parameters where $r \leq 10, i+k \leq 100, k\leq 27$. or $i+k\leq29$ all $r$. A formalism of Wanless as systematized by Pernici is central to this effort.
\end{abstract}

\section{Introduction}
We deal with $r$-regular bipartite graphs with $v = 2n$ vertices.
We let $m_i$ be the number of $i$-matchings.
In \cite{BFP}, Butera, Pernici, and I introduced the quantity $d(i)$, in eq.\
$(10)$ therein,
\begin{equation}\label{eq:1.1}
  d(i) \equiv \ln\biggl( \frac{m_i}{r^i} \biggr)
  - \ln\biggl( \frac{\overline{m}_i}{(v-1)^i} \biggr)
\end{equation}
where $\overline{m}_i$ is the number of $i$-matchings for the complete (not
bipartite complete) graph on the same vertices,
\begin{equation}\label{eq:1.2}
  \overline{m}_i = \frac{v!}{(v-2i)!\,i!\,2^i}
\end{equation}
We here have changed some of the notation from \cite{BFP} to agree with notation
in \cite{Per}.
We then considered $\Delta^k d(i)$ where $\Delta$ is the finite difference
operator, so
\begin{equation}\label{eq:1.3}
  \Delta d(i) = d(i+1) - d(i)
\end{equation}
A graph was defined to satisfy \textit{graph positivity} if all the meaningful
$\Delta^k d(i)$ were non-negative.
That is
\begin{equation}\label{eq:1.4}
  \Delta^k d(i) \ge 0
\end{equation}
for $k = 0,\ldots,v$ and $i = 0,\ldots,v-k$.
We made the conjecture, the 'graph positivity conjecture', supported by some computer evidence,
\begin{conj*}
  As $n$ goes to infinity the fraction of graphs that satisfy graph positivity
  approaches one.
\end{conj*}

We note some of the impressive results of the numerical study of graph positivity in \cite{BFP}.

\begin{enumerate}
	\item All graphs with $v < 14$ satisfy graph positivity.
	\item When $r = 4$, the first violations occur when $v=22$ in 2 graphs of the $2806490$ graphs with $v=22$.
	\item For $r=3$, the fraction of graphs not satisfying graph positivity continuously decreases between $v = 14$ and $v=30$. There is a single violation at $v=14$.
\end{enumerate}

In this paper we study a weaker conjecture than the graph positivity conjecture. We work with ``weak graph positivity''.
\begin{conj*}
 For each $i$ and $k$ one has 
 $$ \Prob(\Delta^kd(i) \geq 0) \underset{n \to \infty}{ \to 1}.$$ 
\end{conj*}
In fact what we prove is:
\begin{thm}\label{thm:1.1}
If $r \leq 10, i+k \leq 100, k\leq 27$, or $i+k\leq29$ all $r$, then 
 $$ \Prob(\Delta^kd(i) \geq 0) \underset{n \to \infty}{ \to 1}.$$ 
\end{thm}
The paper relies heavily on the work of Wanless, \cite{Wan}, and Pernici, \cite{Per}, that gives a nice representation of the $m_i$. The restriction to bipartite graphs is mainly because this restriction is made in \cite{BFP}. In this paper the bipartite nature appears in two places. First, the number of vertices is assumed to be even and second, in eq. (3.8) the lower limit 4 is replaced by 3 if one does not assume the graph is bipartite.

In a previous paper, \cite{6}, we gave a ``near proof" of weak positivity. We assumed the truth of eq. (3.7) and the ``awesome conjecture", Appendix B. Herein we replace the first assumption by Theorem 3.2, and the second by Theorem B.1. Both Theorem 3.2 and Theorem B.1 were proved, in part, using integer valued computer computation, Theorem 3.2 by Pernici, \cite{Per}. Mostly, we follow \cite{6} closely, but with notational changes and organizational rearrangements. 

Clearly, further computer work could likely extend the space of 
parameters for which we can prove weak graph positivity. But the goal of 
future work is to prove eq. (3.7) and the ``awesome conjecture''. Then, maybe to prove graph positivity!

\section{Idea of the Proof}\label{sec:2}
Suppose we want $\Prob(x > y)$ to be large.
We have
\begin{equation}\label{eq:2.1}
  \Prob(x < y) = \Prob(e^x < e^y).
\end{equation}
Set
\begin{equation}\label{eq:2.2}
 (e^x - e^y) \equiv \alpha_0
\end{equation}
and
\begin{equation}\label{eq:2.3}
  \EE(e^x - e^y) \equiv\alpha.
\end{equation}
We will want $ \alpha$ to be positive, and in the present work proving this
positivity will be a major component.
Let
\begin{equation}\label{eq:2.3a}
  \EE\bigl( (e^x - e^y)^2\bigr) - \alpha^2 \equiv \beta.
\end{equation}
Then, assuming $\alpha > 0$,
\begin{equation}\label{eq:2.4}
  \EE\bigl( (e^x - e^y - \alpha)^2\bigr) \ge
   \alpha^2 \, \PP(e^x - e^y < 0).
\end{equation}
And so
\begin{equation}\label{eq:2.5}
  \Prob(e^x < e^y) \le \frac{\beta}{\alpha^2}.
\end{equation}
In our problem $\beta$ and $\alpha$ will be functions of $n$, and we'll want
probability to go to zero with $n$ as $n$ goes to infinity.
\par We turn to the object of study, and perform some simple manipulations,
working from eq.\ \eqref{eq:1.1}
\begin{align}
  &\Prob\bigl( \Delta^k d(i) > 0 \bigr) = \Prob\Biggl( (-1)^k
  \sum_{\ell=0}^k (-1)^\ell \binom{k}{\ell} \, \ln \biggl(
  \frac{m_{i+\ell}}{\overline{m}_{i+\ell}} \cdot \frac{(v-1)^{i+\ell}}{r^{i+\ell}}
  \biggr) > 0 \Biggr)\label{eq:2.6}\\
  &\qquad= \Prob\Biggl( \sum_{\ell \in \mathcal{L}^+} \binom{k}{\ell} \,
  \ln \biggl( \frac{m_{i+\ell}}{\overline{m}_{i+\ell}} \cdot
  \frac{(v-1)^{i+\ell}}{r^{i+\ell}} \biggr) > \sum_{\ell \in \mathcal{L}^-}
  \binom{k}{\ell}\,\ln \biggl( \frac{m_{i+\ell}}{\overline{m}_{i+\ell}} \cdot
  \frac{(v-1)^{i+\ell}}{r^{i+\ell}} \biggr) \Biggr)\label{eq:2.7}
\end{align}
where $\mathcal{L}^+$ is the set of odd $\ell$, $0 \le \ell \le k$, if $k$ is
odd and is the set of even $\ell$, $0 \le \ell \le k$, if $k$ is even, and
$\mathcal{L}^-$ is defined vice versa.
\par Returning to the language of \eqref{eq:2.1}--\eqref{eq:2.5}, we set
\begin{align}
  x &\equiv \sum_{\ell \in \mathcal{L}^+} \binom{k}{\ell}\,\ln \biggl(
  \frac{m_{i+\ell}}{\overline{m}_{i+\ell}} \cdot
  \frac{(v-1)^{i+\ell}}{r^{i+\ell}} \biggr)\label{eq:2.8}\\
  y &\equiv \sum_{\ell \in \mathcal{L}^-} \binom{k}{\ell}\,\ln \biggl(
  \frac{m_{i+\ell}}{\overline{m}_{i+\ell}} \cdot
  \frac{(v-1)^{i+\ell}}{r^{i+\ell}} \biggr)\label{eq:2.9}
  \intertext{and so}
  e^x &= \prod_{\ell \in \mathcal{L}^+} \biggl(
  \frac{m_{i+\ell}}{\overline{m}_{i+\ell}} \cdot
  \frac{(v-1)^{i+\ell}}{r^{i+\ell}} \biggr)^{\binom{k}{\ell}}\label{eq:2.10}\\
  e^y &= \prod_{\ell \in \mathcal{L}^-} \biggl(
  \frac{m_{i+\ell}}{\overline{m}_{i+\ell}} \cdot
  \frac{(v-1)^{i+\ell}}{r^{i+\ell}} \biggr)^{\binom{k}{\ell}}\label{eq:2.11}
\end{align}

For each $k$ we will use these expression to study the large $n$ behavior of $\alpha,\beta,$ and $\beta/\alpha^2$, see (2.3)-(2.6). The proof of Theorem 1.1 follows if $\alpha \geq 0$ and $\beta/\alpha^2 \to 0$. 

Note: throughout often the restriction ``for $n$ large enough'' is understood.

\section{The Work of Wanless and Pernici}\label{sec:3}
In \cite{Wan} Wanless developed a formalism to compute the $m_i$ of any regular
graph.
We here only give a flavor of this formalism, but present some of the consequences we
will use in this paper.
For each $i$ there are defined a set of graphs
$g_{i1},g_{i2},\ldots,g_{i n(i)}$.
Given a regular graph $g$, one computes for each $j$ the number of subgraphs of
$g$ isomorphic to $g_{ij}$, call this $g \sslash g_{ij}$.
Then $m_i$ for $g$ is determined by the $n(i)$ values of $g \sslash g_{ij}$.
We define $M_i$ to be the value of $m_i$ assigned to any graph with all $n(i)$
values of $g \sslash g_{ij}$ zero. Such graphs will exist only for large enough n.
Initially $M_i$ is defined only for such n. But it may be extended as a finite
polynomial in $\frac{1}{n}$ to all non-zero $n$.
$M_i = M_i(r,n)$ is an important object of study to us.
\par In \cite{Per} Pernici systematized the results of Wanless. We now 
present the very non-trivial computational construction of $M_i$ from \cite{Per} and \cite{Wan}. One first defines quantities $u_s(r), s \geq 2$ by 
\begin{equation}
	T_r = \frac{2(r-1)}{2(r-1)-r+r\sqrt{1-4x(r-1)}}. \label{eq:3.1} 
\end{equation}
\begin{equation}
	u_s(r) = [x^s]T_r  \label{eq:3.2}.
\end{equation}
The notation is slightly changed from \cite{Per}. The expression $[x^s]f$ for a series, $f$, in $x$ is defined as the coefficient of $x^s$ in the series $f$. Then one has
\begin{align}
	M_j = [x^j] \exp{\left(nrx - \sum_{s \geq 2}\frac{nu_s(r)}{s} (-x)^s\right)}\label{eq:3.3}.
\end{align}

Defining quantities $a_h(r,j)$ one has the expressions 
\begin{align}
	M_j &= \frac{n^jr^j}{j!}(1+H_j) \label{eq:3.4}\\
	H_j &= \sum_{h = 1}^{j-1} \frac{a_h(r,j)}{n^h}\label{eq:3.5}
\end{align}
that exhibit the structure of $M_j$ especially in so far as powers of $n$. It is important to keep in mind often that 
\begin{align*}
	a_0 &= 1 \\
         a_h &= 0 \text{ if } h \geq j.
\end{align*}

In \cite{Per}, Pernici via a clever formal computation (not rigorous) derives the equations
\begin{align}
   [j^k n^{-h}]\,\ln \biggl( 1 + H_j \biggr) =
    [j^k n^{-h}]\,\ln \biggl( 1 + \sum_{s=1}^{j-1} \frac{a_s(r,j)}{n^s} \biggr) &=
  0, \qquad k \ge h+2 \label{eq:3.6}\\
  [j^{h+1} n^{-h}]\,\ln \biggl( 1 + H_j \biggr)=  [j^{h+1} n^{-h}]\,\ln \biggl( 1 + \sum_{s=1}^{j-1} \frac{a_s(r,j)}{n^s} \biggr)
  &= \frac{1}{(h+1)h} \biggl( \frac{1}{r^h} - 2 \biggr) \label{eq:3.7}
\end{align}
equations (16) and (17) of \cite{Per}. Here $[j^an^{-b}]f$ is an obvious generalization of $[x^s]f$. The status of these equations (\ref{eq:3.6}),(\ref{eq:3.7}) is as follows. First, (\ref{eq:3.6}) is true and in fact we have a stronger result.

\begin{thm}\label{thm:3.1}
	Equation (\ref{eq:3.6}) holds if $M_j$ is calculated from (\ref{eq:3.3}) with any values for the $u_s$, $s\geq 2$, not necessarily the values given by (\ref{eq:3.1}),(\ref{eq:3.2}).
\end{thm}

This is proven in Appendix A, using a major contribution from Robin Chapman, \cite{5}. As to equation (\ref{eq:3.7}) we certainly believe it always holds. But we only claim the cases proved in the next theorem.
\begin{thm}\label{thm:3.2}
	For $r \leq 10$ and $j\leq100$ , and for $j\leq29$ and all $r$, equation (\ref{eq:3.7}) holds.
\end{thm}

Pernici proved this for $r \leq 10$ and $j\leq100$ using computer computation. See the end of Section 2 in \cite{Per}. The additional region was proved by me in [10], also by computer computation. We note that Theorem \ref{thm:3.2} yields inductively values of $a_k(r,j)$ for $k \leq 100$, $j\leq 100$, $r\leq 10$. One should see \cite{Per} for a logical development of this and Theorem \ref{thm:3.2}. A useful inductive deduction from Theorem 
\ref{thm:3.1} is that $a_k$ is a polynomial of degree at most $2k$.

\par We set $M_0 = 1$ and $M_s = 0$ if $s < 0$.
Then $m_j$ is recovered from $M_j$ by the formula
\begin{equation}\label{eq:3.8}
  m_j = \exp \biggl( \sum_{s \ge 4} \frac{\varepsilon_s}{2s} (-\hat{x})^s \biggr)
  M_j
\end{equation}
eq.\ $(11)$ of \cite{Per}. Here
\begin{equation}\label{eq:3.9}
  \hat{x} M_j = M_{j-1}
\end{equation}
$\varepsilon_s$ for a graph $g$ is a linear function of a finite number of $g
\sslash \ell_i$, $\ell_i$ a set of given graphs, the 'contributors'.
The only thing we need to know is that for any given product of $\varepsilon_s$'s, $\prod_i
\varepsilon_{s(i)}$, one has that
\begin{equation}\label{eq:3.10}
  \EE\biggl( \prod_i \varepsilon_{s(i)} \biggr) \le C
\end{equation}
i.e.\ it is a bounded function of $n$. Here as everywhere in this paper the expectation
is the average value of the function
over all r regular bipartite graphs of order $2n$.

The result one needs to see this is that the number of $s$-cycles are
independent Poisson random variables of finite means in the fixed $r$, $n$ goes
to infinity limit, \cite{Bol}. The needed extension to the bipartite case was done
by Wormald, [9].
One then uses the fact that the $\ell_i$ and $g_{ij}$ graphs discussed above all
are either single cycle or multicycle in nature.

Working from eq.\ $(3.8)$ one can arrange the resultant terms arising 
into the following expression for $m_j$
\begin{align}
  m_j &= \frac{n^jr^j}{j!} (1+\hat{H_j)} \label{eq:3.11}\\
 \hat {H_j }&= \sum_{h=1}^{j-1} \frac{\hat{a_h}(r,j,\{\epsilon_i\})}{n^h} \label{eq:3.12}
\end{align}
$m_j$ is a function on graphs, eq.\ $(3.8)$ or eq.\ $(3.11)$-$(3.12)$ in turn expresses $m_j$ as a polynomial
in the $\{\epsilon_i\}$, these also functions on the graphs. We will be dealing with expectations of
polynomials in the $\{m_j\}$, for example eq.\ $(4.3)$. We make the important observation that,
for the sum in eq.\ $(3.12)$ appearing in an expectation,
the $n$ dependence of the $\epsilon_i$ does not effect the formal expected asymptotic expansion
by powers of $1/n$, from the discussion surrounding
eq.\ $(3.10)$.

\section{Some simple reorganization}
We define
\begin{equation}\label{eq:4.1}
  1 + K_i \equiv \frac{(v-1)^i}{r^i} \cdot \frac{(v-2i)!\,i!\,2^i}{v!} \cdot
  \frac{r^in^i}{i!}
\end{equation}
using notably eq.\ \eqref{eq:1.2}.
Then with
\begin{equation}\label{eq:4.2}
  \alpha_0=\Biggl(
    \prod_{\ell \in \mathcal{L}^+} \bigl(
      (1+\hat{H_{i+\ell}})(1+K_{i+\ell})
    \bigr)^{\binom{k}{\ell}}
    - \prod_{\ell \in \mathcal{L}^-} \bigl(
      (1+\hat{H_{i+\ell}})(1+K_{i+\ell})
    \bigr)^{\binom{k}{\ell}}
  \Biggr)
\end{equation}
$\alpha$ becomes
\begin{equation}\label{eq:4.3}
  \alpha=\EE(\alpha_0)
 \end{equation}

Further we set
\begin{equation}\label{eq:4.4}
  1 + K_i \equiv e^{G_i}
\end{equation}
where
\begin{align}
  G_i &\equiv G_{i,1} + G_{i,2} + G_{i,3} + G_{i,4} + G_{i,5} \label{eq:4.5}\\
  G_{i,1} &\equiv i\,\ln\Bigl( 1 - \frac{1}{2n} \Bigr) \label{eq:4.6}\\
  G_{i,2} &\equiv (2n-2i)\,\ln\Bigl( 1 - \frac{i}{n} \Bigr) \label{eq:4.7}\\
  G_{i,3} &\equiv 2i \label{eq:4.8}\\
  G_{i,4} &\equiv \frac{1}{2}\,\ln\Bigl( 1 - \frac{i}{n} \Bigr) \label{eq:4.9}\\
  G_{i,5} &\equiv \sum_{j\ \text{odd}} c_j \biggl( \frac{1}{n^j} -
  \frac{1}{(n-i)^j} \biggr) \label{eq:4.10}
\end{align}
We have used the Stirling series to expand $\ln n!$.
We also note that for example $c_1 = -\frac{1}{24}$.
 $K_i$ is easily developed as a series in inverse powers of $n$.
 
 \textit{Valuable Observation} The convergence problem for series, except inside
 expectation values, is trivial, since one deals with $r$, $j$, and $\epsilon_i$ ( taken as a
 number ) fixed and $n$ large enough. BUT, the only expectations we take are of $\alpha_0$
 and $\alpha_0^{2}$ ( for $\beta$ ). And, see (4.2) and the discussion after (3.11)-(3.12), these
 both are finite polynomials in the $\{\epsilon_i\}$! ( So to study $\alpha_0$ and $\alpha$
 it is a good idea to expand $\alpha_0$ in the formal series
 in powers of $\frac{1}{n}$ taking the coefficients of the terms through $\frac{1}{n^{k-1}}$ from
 eq.(6.2) and the rest of the terms from eq.(4.2).)

\section{First Identity}\label{sec:5}
We start with some simple manipulations
\begin{equation}\label{eq:5.1}
  \prod_i (1+x_i)^{e_i} = e^{\sum e_i\ln(1+x_i)} = 1 + \Bigl( \sum e_i\ln(1+x_i)
  \Bigr) + \frac{1}{2!} \Bigl( \sum e_i\ln(1+x_i) \Bigr)^2 + \cdots
\end{equation}
With the notation
\begin{equation}\label{eq:5.2}
  (1 + \hat{H_i)}(1 + K_i) \equiv 1 + U_i
\end{equation}
we substitute $U_i$ for $x_i$ and $\binom{k}{\ell}$ for $e_i$ in \eqref{eq:5.1}
\begin{equation}\label{eq:5.3}
  \prod_{\ell \in \mathcal{L}^+} (1 + U_{i+\ell})^{\binom{k}{\ell}} = 1 + t _+ +
  \frac{1}{2} t_+^2 + \frac{1}{3!} t_+^3 \cdots
\end{equation}
where
\begin{equation}\label{eq:5.4}
  t _+\equiv \sum_{\ell \in \mathcal{L}^+} \binom{k}{\ell} \biggl( U_{i+\ell} -
  \frac{1}{2} (U_{i+\ell})^2 + \frac{1}{3} (U_{i+\ell})^3 \cdots \biggr)
\end{equation}
The First Identity consists of \eqref{eq:5.3} and \eqref{eq:5.4} and the same
expressions with $\mathcal{L}^+$, $t_+$ replaced by $\mathcal{L}^-$, $t_-$.

\section{Second Identity }\label{sec:6}
Throughout this section we understand the restrictions
\begin{equation}
	r\leq 10, k\geq 2, k\leq 27, i+k \leq 100 \label{eq:6.1}
\end{equation}
For convenience we introduce
\begin{align}
	F_i \equiv (1+\hat{H_i})(1+K_i) \quad  \label{eq:6.2}.
\end{align}
\begin{thm}[Second Identity]\label{thm:6.1}
	For $r\leq 10, k\geq 2, k\leq 27, i+k \leq 100$ 
	\begin{equation}\label{eq:6.3}
    \sum_{\ell=0}^k \binom{k}{\ell} (-1)^{\ell+k} \biggl[ \frac{1}{n^{k-1}}
    \biggr] \sum_{m=1}^{k-1} (-1)^{m+1} \frac{1}{m} \bigl( F_{(i+\ell)} - 1
    \bigr)^m = \frac{(k-2)!}{r^{k-1}}
  \end{equation}
\end{thm}
\begin{thm}\label{thm:6.2}
	For $r\leq 10, k\geq 2, k\leq 27, i+k \leq 100$ 

	\begin{align}\label{eq:6.4}
		\left[ \frac{1}{n^{k-1}}\right]\ln(F_i)
	\end{align}
	has highest power of $i = i^k$, and this term is
	\begin{align}\label{eq:6.5}
	\frac{(k-2)!}{k!}\frac{i^k}{r^{k-1}}
	\end{align}
\end{thm}
For example for $k = 3$
\begin{equation}\label{eq:6.6}
  \biggl[\frac{1}{n^2}\biggr] \ln(F_s) = -\frac{1}{12} \frac{s\bigl(3r^2s - 3r^2
  - 12rs - 2s^2 + 12r + 9s - 7\bigr)}{r^2}
\end{equation}
\par We now note
\begin{equation}\label{eq:6.7}
  \sum_{\ell=0}^k \binom{k}{\ell} (-1)^{\ell+k} \ell^d = \begin{cases}
    \hfil 0 & d < k\\
    k! & d = k
  \end{cases}
\end{equation}
that follows from
\begin{equation}\label{eq:6.8}
  \sum_{\ell=0}^k \binom{k}{\ell} (-1)^{\ell+k} \ell^d = \Delta^k i^d
\end{equation}
which like $(\frac{d}{dx})^k x^d$ has values in \eqref{eq:6.7}.
From eq.\ \eqref{eq:6.7} one can deduce that Theorem \ref{thm:6.1} follows from
Theorem \ref{thm:6.2}, which we proceed to prove.
\begin{equation}\label{eq:6.9}
  \biggl[ \frac{1}{n^{k-1}} \biggr] \ln(F_i) = 
  \biggl[ \frac{1}{n^{k-1}} \biggr] \ln(1+\hat{H_i})
  + \biggl[ \frac{1}{n^{k-1}} \biggr] \ln(1+K_i)
\end{equation}
From Theorem B.1 with $k \ge 2$ we see that the highest
power of $i$ in $[\frac{1}{n^{k-1}}]\ln(1+\hat{H_i)}$ is $k$ and its coefficient is
\begin{equation}\label{eq:6.10}
  \frac{1}{k(k-1)} \biggl( \frac{1}{r^{k-1}} - 2 \biggr)
\end{equation}
To study $[\frac{1}{n^{k-1}}]\ln(1+K_i)$ we turn to equations
(4.4)-(4.10).
We note the highest power of $i$ arises from the expansion of the term
$G_{i,2}$, eq.(4.7), and $[\frac{1}{n^{k-1}}]\ln(1+K_i)$ has highest
power of $i$ equal $k$ and its coefficient is
\begin{equation*}
  \frac{2}{k(k-1)}
\end{equation*}
So
\begin{equation}\label{eq:6.11}
  \biggl[ \frac{1}{n^{k-1}} \biggr] \ln(F_i) = \frac{(k-2)!}{k!}
  \frac{1}{r^{k-1}} i^k
\end{equation}
Quod erat demonstrandum.

\section{$k=1$ and $k=0$}\label{sec:7}
Not only is $k=1$ the first case, but it is different from $k\geq 2$ in some essential ways. We proceed to compute $\alpha$ for $k=1$. From eq. \eqref{eq:4.2} we have 
\begin{align}\label{eq:7.1}
	\alpha_0 = (1+\hat{H_{i+1}})(1+K_{i+1})-(1+\hat{H_i})(1+K_i)
\end{align}
Following the development in \cite{Per}, from eq. (17) through the end of Section 2, using Theorem B.1 of this paper, one gets 
\begin{equation} \label{eq:7.2}
	\hat{H_i} \sim H_i \quad r\leq 10, i\leq 100
\end{equation}
i.e. they have the same leading term. 

Taking $H_i$ in this range, from eq(18) and eq(45) of \cite{Per}
\begin{equation}\label{eq:7.3}
	\hat{H_i} = i(i-1)\left(-1+\frac{1}{2r}\right)\frac{1}{n} + \mathcal{O}\left(\frac{1}{n^2}\right) \quad r\leq 10, i\leq 100.
\end{equation}

In such asymptotic series bounds we treat the $\varepsilon_i$ as constants. 

Using eq. \eqref{eq:4.4}-\eqref{eq:4.10} one gets
\begin{align}\label{eq:7.4}
	K_i = (i^2-1)\frac{1}{n} + \mathcal{O}\left(\frac{1}{n^2}\right).
\end{align}
There follows

\begin{thm}\label{thm:7.1}
	For $k=1$
	\begin{align}\label{eq:7.5}
		\alpha_0 = \frac{i}{rn}+\mathcal{O}\left(\frac{1}{n^2}\right) \quad r\leq 10, i\leq 99
	\end{align}
\end{thm}

One easily gets

\begin{thm}\label{thm:7.2}
	For $k=0$
	\begin{align}\label{eq:7.6}
		\alpha_0 = 1+\mathcal{O}\left(\frac{1}{n}\right) \quad r\leq 10, i\leq 100
	\end{align}

\end{thm}

\section{$k\geq 2$}\label{sec:8}
Throughout this section we enforce the restrictions of eq. \eqref{eq:6.1}
\begin{align*}
	r\leq 10, k\geq 2, i\geq 0, i+k \leq 100, k\leq 27.
\end{align*}

The goal of this section is proving the following theorem.
\begin{thm}\label{thm:8.1}
  For $k \ge 2, r\leq 10, i+k \leq 100, k \leq 27$
  \begin{equation}\label{eq:8.1}
 \alpha_0 = \frac{(k-2)!}{r^{k-1}n^{k-1}} + \OO(\frac{1}{n^{k}})
  \end{equation}
  
\end{thm}
From Section \ref{sec:5} using the First Identity we have
\begin{equation}\label{eq:8.2}
  \alpha_0 = \biggl( \Bigl( 1 + t_+ + \frac{1}{2}t_+^2 + \cdots\Bigr) -
  \Bigl( 1 + t_- + \frac{1}{2}t_-^2 + \cdots\Bigr) \biggr)
\end{equation}
where with
\begin{equation}\label{eq:8.3}
  1 + U_i = (1 + \hat{H_i})(1 + K_i)\\
\end{equation}
one defines
\begin{equation}
 t_+ = \sum_{\ell \in \mathcal{L}^+} \binom{k}{\ell}
  \Bigl( U_{i+\ell} - \frac{1}{2} (U_{i+\ell})^2 + \frac{1}{3} (U_{i+\ell})^3
  \cdots \Bigr)\label{eq:8.4}\\
\end{equation}
\begin{equation}
  t_- = \sum_{\ell \in \mathcal{L}^-} \binom{k}{\ell}
  \Bigl( U_{i+\ell} - \frac{1}{2} (U_{i+\ell})^2 + \frac{1}{3} (U_{i+\ell})^3
  \cdots \Bigr)\label{eq:8.5}
\end{equation}
We treat the terms written explicitly in \eqref{eq:8.2}; first the linear terms.
\begin{align}
  \biggl[ \frac{1}{n^d} \biggr] (t_+ - t_-) &= \biggl[ \frac{1}{n^d} \biggr]
  \sum_\ell \binom{k}{\ell} (-1)^{k+\ell} \Bigl( U_{i+\ell} - \frac{1}{2}
  (U_{i+\ell})^2 + \frac{1}{3} (U_{i+\ell})^3 \cdots \Bigr)\label{eq:8.6}\\
  &= \begin{cases}
    \hfil 0 & d < k - 1\\
    \displaystyle\frac{(k-2)!}{r^{k-1}} & d = k-1
  \end{cases}\label{eq:8.7}
\end{align}
by the Second Identity, Theorem \ref{thm:6.1}. In applying $\bigg[\frac{1}{n^{d}}\bigg]$
we treat the $\epsilon_i$ as constants.
Next we want to prove that the higher powers of $t$'s make no contribution in
\eqref{eq:8.2}!
This is amazing when one first sees it. Explicitly dealing only with the quadratic terms,
\par We want to show
\begin{equation}\label{eq:8.8}
  \biggl[ \frac{1}{n^d} \biggr] (t_+^2 - t_-^2) = 0 \quad \text{for} \quad d
  \le k - 1
\end{equation}
We proceed by looking at the powers of $\frac{1}{n}$.
\begin{equation}\label{eq:8.9}
  \biggl[ \frac{1}{n^d} \biggr] (t_+^2 - t_-^2) = \sum_{s=1}^{d-1} \Biggl[
    \biggl( \biggl[ \frac{1}{n^s} \biggr] t_+ \biggr)
    \biggl( \biggl[ \frac{1}{n^{d-s}} \biggr] t_+ \biggr)
  - \biggl( \biggl[ \frac{1}{n^s} \biggr] t_- \biggr)
    \biggl( \biggl[ \frac{1}{n^{d-s}} \biggr] t_- \biggr) \Biggr]
\end{equation}
All we need to complete a proof of \eqref{eq:8.8} is to show
\begin{equation}\label{eq:8.10}
  \biggl[ \frac{1}{n^s} \biggr] t_+ = \biggl[ \frac{1}{n^s} \biggr] t_- \qquad
  1 \le s \le d - 1
\end{equation}
But this follows from \eqref{eq:8.6},\eqref{eq:8.7} above.
Pretty neat.

\section{Completion}\label{sec:9}
In this section, we always have $r \leq 10, i+k \leq 100, k \leq 27$. The information we need from the 
calculations of this paper are Theorem \ref{thm:7.1}, Theorem \ref{thm:7.2}, and Theorem \ref{thm:8.1}. From
these respectively, we get:
\begin{enumerate}
	\item For $k=1,i\neq 0$:
		\begin{align}\label{eq:9.1}
			\alpha&\geq \frac{c_1}{n} \quad c_1 \text{ positive}\\ 
		\label{eq:9.2}
			\beta &\leq \frac{c_2}{n^4} 
		\end{align}
	\item For $k=0$:
		\begin{align}
			\alpha&\geq c_3 \quad \text{$c_3$ positive} \\
			\beta &\leq \frac{c_4}{n^2}
		\end{align}

	\item For $k \geq 2$:
		\begin{align}
			\alpha &\geq \frac{c_{5,k}}{n^{k-1}} \quad \text{$c_{5,k}$ positive} \\
			\beta &\leq \frac{c_{6,k}}{n^{2k}}
		\end{align}
\end{enumerate} Referring to Section \ref{sec:2}, Theorem \ref{thm:1.1} follows from the fact that $\beta/\alpha^2$ goes to zero as $n$ goes to infinity in each case. The key observation in getting the estimate for $\beta$ is to notice that in 
\begin{equation}\label{eq:9.7}
	\beta = E(\alpha_0^2) - \left(E(\alpha_0)\right)^2
\end{equation}
the lead term and the second term in an expansion in powers
of $\frac{1}{n}$ for the two terms on the right side of \eqref{eq:9.7} cancel since the lead term in the expansion of $\alpha_0$ is independent of the $\varepsilon_i$!

\appendix
\section{Appendix}\label{app:A}
The organization of this Appendix is as follows. Section A.1 presents Theorem A.1, which is actually a restatement of Theorem \ref{thm:3.1}. Section A.2 presents Theorem A.2, a theorem about Stirling numbers of the first king,\cite{5}, \cite{7}.
In Section A.3 we prove the equivalence of Theorem A.1 and Theorem A.2,\cite{5}, \cite{7}. Thus a proof of Theorem \ref{thm:3.1}.

Historically, I wanted to prove Theorem A.1. I found that Theorem A.1 implied Theorem A.2. I put Theorem A.2 as a conjecture, onto the web, \cite{7}. Robin Chapman saw this paper and proved the theorem, \cite{5}, by a clever little argument.

\subsection*{Section A.1}
Let $\tilde{u_s}, s\geq 2$ be arbitrary numbers. Define
\begin{equation}\label{eq:a.1}
	\tilde{M_j}(n) \equiv [x^j]\exp\left(nx+\sum_{s\geq 2} n \tilde{u_s} x^s\right)
\end{equation}

\begin{equation}\label{eq:a.2}
	\tilde{M_j} \equiv \frac{n^j}{j!} \sum_{n=0}^{j-1} \frac{\tilde{a}_n(j)}{n^n}
\end{equation}
Then we have
\begin{thm}\label{thm:a.1}
	\begin{equation}\label{eq:a.3}
	[j^kn^{-h}] \ln\left(1+\sum_{s=1}^{j-1}\frac{\tilde{a}_s}{n^s}\right) = 0 \quad k\geq h+2\end{equation}
\end{thm}	
\subsection*{Section A.2}
The (unsigned) Stirling numbers of the first kind, ${\stirlingi{a}{b}}$, are defined by
\begin{equation}
	x(x+1)\dots(x+n-1) = \sum_{k=1}^n \stirlingi{n}{k} x^k 
\end{equation}
It is easy to show $\stirlingi{n}{n-w}$ is a polynomial in $n$ of degree $2w$. So we may naturally define
$\stirlingi{x}{x-w}$ for any number $x$ by extending the
domain of the polynomial. We set
\begin{align}
	P_w(x) = \stirlingi{x}{x-w}
\end{align}

Now we give ourself an integer $g\geq 2$, an integer $w$, $0 \leq w \leq g-2$ and a set of $g$ distinct numbers,
\begin{align}
	S = \{c_1,\dots,c_g\}.
\end{align}
We define a \underline{configuration} as a sequence of non-empty subsets of $S$:
\begin{align}
	S_1,\dots,S_m
\end{align}
that are disjoint with union $S$,i.e.
\begin{align}
	S_i \neq \emptyset, S_i \cap S_j = \emptyset \text{ if } i\neq j \quad \text{ and } \bigcup_{i=1}^m S_i = S.
\end{align}
For a configuration we define
\begin{align}
	t_i = \sum_{c_j \in S_i} c_j \quad i=1,\dots,m
\end{align}
A \underline{weighted configuration} is a configuration as above for which each $S_i$ is assigned a non-negative integer $w_i$, its weight, with the restriction
\begin{align}
	\sum_{i=1}^m w_i = w.
\end{align}
Such a weighted configuration has an \underline{evaluation} defined as
\begin{align}
	(-1)^m \frac{1}{m} \prod_{i=1}^m P_{w_i}(t_i)
\end{align}

\begin{thm}\label{thm:a.2}
The sum over all distinct weighted configurations of their evaluations is zero.
\end{thm}

\subsection*{Section A.3}
\begin{thm}
Theorem A.1 implies Theorem A.2
\end{thm}

Given integer $g \geq 2$ and integer $w, 0 \leq w \leq g-2$ and a set of $g$ distinct numbers 
\begin{align}
	S = \{c_1,\dots,c_g\}
\end{align}
we want that the associated set of weighted configurations have a sum of their evaluations equal zero.
Of course, read Section A.2. 
We first prove this in the case the set $S$ has all the $c_i$ positive integers, each $\geq 2$. 

We set 
\begin{align}
	k &= \sum c_i - w \label{eq:a.13}\\
	\intertext{and}
	h &= \sum c_i - g\label{eq:a.14}
\end{align}
so that $k-h = g-w \geq 2$. We set all the $\tilde{u}_s$ to zero in $\eqref{eq:a.1}$, 
except for values of $s$ equal to one of the $c_i$. Then the left side of \eqref{eq:a.3} is identically zero as a function of the $\tilde{u}_s$. The coefficient of the term $\tilde{u}_{c_1}\tilde{u}_{c_2}\dots\tilde{u}_{c_g}$ in the left side of \eqref{eq:a.3} $=0$ is just the statement of Theorem \ref{thm:a.2}.

Since the sum over evaluations, considered in Theorem \ref{thm:a.2} is a polynomial in the $c_i$ of degree
$2w$, knowing it is true for all integer $c_i$, distinct,
$\geq 2$, implies that it is true for all $c_i$. 

\begin{thm}
	Theorem \ref{thm:a.2} implies Theorem \ref{thm:a.1}
\end{thm}

We introduce the notation $\Phi(g,w,\{c_i\})$ for the sum in Theorem \ref{thm:a.2}. Then with $c_1,\dots,c_g$ distinct integers $\geq 2$ what the reader computed in showing Theorem \ref{thm:a.1} implies Theorem \ref{thm:a.2} above was the equivalence of 
\begin{align}
	\Phi(g,w,\{c_i\}) 
\end{align}
with
\begin{align}	[\tilde{u}_{c_1}\tilde{u}_{c_2}\dots\tilde{u}_{c_g}][j^kn^{-h}] \ln\left(1+\sum_{s=1}^{j-1}\frac{\tilde{a}_s}{n^s}\right)
\end{align}
The equations (\ref{eq:a.1},\ref{eq:a.2},\ref{eq:a.13},\ref{eq:a.14}) being understood. Thus we know that in \eqref{eq:a.3} the terms in a development into a power series in the $\tilde{u}_i$ for which there are are no repeated variables (terms linear in each $\tilde{u}_i$) are zero.

We consider a term with repeated indices
\begin{equation}
	[\tilde{u}_{c_1}^{\alpha_1}\tilde{u}_{c_2}^{\alpha_2}\dots\tilde{u}_{c_g}^{\alpha_g}][j^kn^{-h}]\ln\left(1+\sum_{s=1}^{j-1}\frac{\tilde{a}_s}{n^s}\right)
\end{equation}

\begin{align}
	k=\sum \alpha_i c_i - w \\
	h = \sum \alpha_i c_i -g 
\end{align}
where at least one of the $\alpha_i$ is $>1$. We now take advantage of the fact that in Theorem A.2 the $c_i$
do not have to be integers. We consider

\begin{align}
	\Phi\left(\sum \alpha_i,w,\{c_{ij}\}\right)
\end{align}
where for each $i$ the range of $j$ is $1,2,\dots,\alpha_i$
the $c_{ij}$ distinct. It is natural to take a limit where each $c_{ij}$ approaches $c_i$ through a limit where 
the $c_{ij}$ are distinct
\begin{align}\label{eq:a.21}
\widetilde{\lim_{\{c_{ij}\to c_i\}}}\Phi\left(\sum\alpha_i,w,\{c_{ij}\}\right).
\end{align}
The tilde over the limit means that the $c_{ij}$ are kept distinct through the limit process. Of course, \eqref{eq:a.21} is zero. 

What is true, proving the implication of this section of the appendix is that 
\begin{align}\label{eq:a.22}
	\widetilde{\lim_{\{c_{ij}\to c_i\}}}\Phi\left(\sum\alpha_i,w,\{c_{ij}\}\right) =\left( \prod_i (\alpha_i)!\right) [\tilde{u}_{c_1}^{\alpha_1}\tilde{u}_{c_2}^{\alpha_2}\dots\tilde{u}_{c_g}^{\alpha_g}][j^kn^{-h}]\ln\left(1+\sum_{s=1}^{j-1}\frac{\tilde{a}_s}{n^s}\right)
\end{align}

One can think that on both sides of the equation one is dealing for each $i$ with $\alpha_i$ ``identical but distinct" integers $c_i$! Whereas the implication Theorem \ref{thm:a.2} $\implies$ Theorem \ref{thm:a.1} required a painfully fussy calculation (but just a calculation), the understanding of \eqref{eq:a.22} and the inverse of the 
implication can be reasoned through in your head, 
no calculation. At the end of the day, that the
$c_i$ in Theorem \ref{thm:a.2} need not be an integer
makes all the difference!

\section{Appendix}
In this appendix, we deal with ``the awesome conjecture", Conjecture B.1. Section B.1 presents this conjecture. Also presented there is the weaker theorem, Theorem B.1, that is a partial substitute for ``the awesome conjecture". The proof of Theorem B.1 is presented in Section B.2 and Section B.3. ``The awesome conjecture"
was first put forth in \cite{6}. Lemma B.1 in Section B.2 proved by a computer, is an example of ``the Genius conjectures" of \cite{8}.

\subsection*{Section B.1}
\begin{conj*}[The Awesome Conjecture]
	Let $z_i \geq 2$ be positive integers. We set:
\begin{equation}\label{eq:b.1}
F =  \sum_{s \ge 0} \frac{a_s(r,j)}{n^s} +\sum _i c_i j (j-1)\cdots (j-z_i+1)\frac{1}{n^{z_i} r^{z_i}} \sum_{s \ge 0} \frac{a_s(r,j-z_i)}{n^s} 
\end{equation}
Then we conjecture:
\begin{align}
  [j^k n^{-h}]\,\ln ( F ) &=
  0, \qquad k \ge h+2 \label{eq:b.2}\\
  [j^{h+1} n^{-h}]\,\ln ( F )
  &= \frac{1}{(h+1)h} \biggl( \frac{1}{r^h} - 2 \biggr) \label{eq:b.3}
\end{align}
\end{conj*}

Compare with eq.(3.6)-(3.7), the $a$'s are the same. 

In the next two section we will prove Theorem B.1 as much as we could easily prove \underline{by computer computation}.

\begin{thm}\label{thm:b.1}
Let $z_i \geq 4 $ be positive integers. We set 
\begin{equation}\label{eq:b.4}
F =  \sum_{s \ge 0} \frac{a_s(r,j)}{n^s} +\sum _i c_i j (j-1)\cdots (j-z_i+1)\frac{1}{n^{z_i} r^{z_i}} \sum_{s \ge 0} \frac{a_s(r,j-z_i)}{n^s} 
\end{equation}
Then if $r \leq 10, h \leq 26, j\leq 100$

\begin{align}
  [j^k n^{-h}]\,\ln ( F ) &=
  0, \qquad k \ge h+2 \label{eq:b.5}\\
  [j^{h+1} n^{-h}]\,\ln ( F )
  &= \frac{1}{(h+1)h} \biggl( \frac{1}{r^h} - 2 \biggr) \label{eq:b.6}
\end{align}
\end{thm}

( If the restriction $z_i \geq 4$ is changed to $ z_i \geq 2$  the upper limits on $k$ in the 
paper are changed from $k \leq 27$ to $k \leq 21$. )

\subsection*{Section B.2}
The content of eq. \eqref{eq:3.3},\eqref{eq:3.4},\eqref{eq:3.5} we are now after is contained in 
\begin{align}
	M_j = [x^j] \exp{\left(nrx - \sum_{s \geq 2}\frac{n u_s}{s} (-x)^s\right)}\label{eq:b.7}
\end{align}
\begin{align}
	\left(1+\sum_{h=1}^{j-1} \frac{a_h(r,j)}{n^h}\right)=\frac{j!}{n^jr^j}M_j.
\end{align}
This we view as a vehicle to pass from a given sequence of functions $u_2,u_3,\dots$ to a sequence of functions $a_1,a_2,\dots$. Each of the functions is a function of $r$ and $j$.
We introduce a function, $\mathcal{F}$, so we may write
\begin{align}
	\{a_h\} = \mathcal{F}(\{u_s\})\label{eq:b.9}.
\end{align}
Also, eq. \eqref{eq:a.1} and \eqref{eq:a.2} are of this nature
\begin{align}
	``\{\tilde{a_h}\} = \mathcal{f}(\{\tilde{u_s}\})''\label{eq:b.10}
\end{align}
where the quotes we have put in \eqref{eq:b.10} are to point out we have to make some simple redefinitions of the $\tilde{a_h}$ and $\tilde{u_s}$ for the equation to hold. 
We are now in a position to state the result of this section.

\begin{lemma}\label{lemma:b.1}
Let $z_i$ be an integer $\geq 2$. Then if 
\begin{equation}\label{eq:b.11}
F =  \sum_{s \ge 0} \frac{a_s'(r,j)}{n^s} +\sum _i c_i j (j-1)\cdots (j-z_i+1)\frac{1}{n^{z_i} r^{z_i}} \sum_{s \ge 0} \frac{a_s'(r,j-z_i)}{n^s} 
\end{equation}
where
\begin{align}
	\{a_s'\} = \mathcal{F}(\{u_s'\})\label{eq:b.12},
\end{align}
there are $xx_s$ such that 
\begin{align}
	F = \sum_{s \geq 0} \frac{a_s''}{n^s}\label{eq:b.13}
\end{align}
where
\begin{align}
	\{a_s''\} = \mathcal{F}(\{u_s'+xx_s\}).\label{eq:b.14}
\end{align}
Moreover, if each $u_s'$ is bounded as $j \to \infty$, and the $z_i \geq 4$ then
\begin{align}
	xx_s = \mathcal{O}(1/j) \text{ if } s \leq 27\label{eq:b.15}.
\end{align}

\end{lemma}
The existence of $xx_s$ such that \eqref{eq:b.14} holds is trivial. The statement after the moreover, \eqref{eq:b.15} is the real content of the theorem. 
\eqref{eq:b.15} holding without the restriction $s \leq 27$ would be a consequence of the ``the Genius Conjectures'' [8]. The present result is a computer computation.

The computation was done in Maple. We 
used less than an hour on a simple desktop computer. 

The outuput of the program included $xx_s$ in the form
\begin{align}
	xx_s = P_s/Q_s \label{eq:b.16}
\end{align}
where $Q_s$ is a polynomial in $j$ and $P_s$ is a polynomial in the $\{u_s\},j,\{c_i\}$. $r$ is irrelevant and may be set equal to one. It also checked that the degree in $j$ of $P_s$ is less
than the degree of $Q_s$, All for $s\leq 27$. Of course, we designed the line of development in Appendix B to make efficient use of the computers computational power.

\subsection*{Section B.3}
We again introduce some functionals. For the basic equations
\begin{align}
	M_j = [x^j] \exp{\left(nrx - \sum_{s \geq 2}\frac{nu_s'}{s} (-x)^s\right)}\label{eq:b.17}
\end{align}
\begin{align}
	(1+H_j') = \frac{j!}{n!r^j} M_j \label{eq:b.18}
\end{align}
we let $\mathcal{H}$ be the functional so determining a map from the $u_s'$ to the $H_j'$
\begin{align}
	(1+H'_j) = \mathcal{H}(\{u_s'\}).\label{eq:b.19}
\end{align}
For the relation
\begin{equation}\label{eq:b.26}
F =  \sum_{s \ge 0} \frac{a_s'(r,j)}{n^s} +\sum _i c_i j (j-1)\cdots (j-z_i+1)\frac{1}{n^{z_i} r^{z_i}} \sum_{s \ge 0} \frac{a_s'(r,j-z_i)}{n^s} 
\end{equation}
we set the functional $\mathcal{G}$ such that 
\begin{align}
	F = \mathcal{G}(\{a_s'\})\label{eq:b.21}
\end{align}
so that with $\mathcal{F}$ from eq. \eqref{eq:b.7}-\eqref{eq:b.9} one has
\begin{align}
	F=\mathcal{G}\circ\mathcal{F}(\{u_s'\}).\label{eq:b.22}
\end{align}
We collect some basic information we want to emphasize.
\begin{enumerate}[I.)]
	\item With $(1+H_j') = \mathcal{H}(\{u_s'\})$ one has that $\ln(1+H_j')$ is a polynomial in $\frac{1}{n},j,$ and the $\{u_s'\}$. So $r$ may appear in the coefficients of the polynomial.
	\item With $F=\mathcal{G}\circ\mathcal{F}(\{u_s'\})$ one has that $(F)$ is a polynomial in $\frac{1}{n},j,$ and the $\{u_s'\}$.
	\item Given a set $\{u_s'\}$ there are $\{xx_s'\}$ such that 
		\begin{align}
 ln  (\mathcal{H}(\{(u_s'+xx_s')\}))= \mathcal{G}\circ\mathcal{F}(\{u_s'\})\label{eq:b.23}
		\end{align}
and if each $u_s'$ is bounded as $j$ goes to infinity,  and the $z_i \geq 4$ then
\begin{align}
	xx_s' = \mathcal{O}(1/j) \text{ if } s \leq 27 \label{eq:b.24}
\end{align}
(Note: $u_s'$ may be a function of $r$ and $j$.)
The last result is from Lemma \ref{lemma:b.1}.

With I, II, and III, the proof of Theorem \ref{thm:b.1} is within mental grasp. (Of course with Theorem \ref{thm:3.1} and Theorem \ref{thm:3.2}.)
\end{enumerate}

\section{Appendix}
It is sad to know that Robin Chapman is now dead. He did not complete preparation of \textit{ A Stirling number identity } of which
he sent me a rough draft. Nor do we know if he would have sought to publish it. The proof we present follows
the lines of his rough draft. This presentation is centered about proving eq.(C.1) below. It was I think particularly clever
of Robin Chapman to realize this was the 'right' equivalent formulation of our Theorem A.2 ( presented as a conjecture
on the web ) to address.
\bigskip

\noindent
1) We note the Stirling Number of the first kind $\stirlingi{n}{k}$ is the number of permutations of $n$ objects with $k$ cycles.And
the Stirling Number of the second kind $\stirlingii{n}{k}$ is the number of partitions of $n$ objects into $k$ disjoint non-empty sets.
\bigskip

\noindent
2) We set $Sym(A)$ to be the set of permutations of set $A$. The "weight" of a permutation $\sigma \in \symm(A)$ as
\begin{equation*}
 w(\sigma) = n - c(\sigma)
\end{equation*}
the number of elements in $A$ minus the number of cycles in $\sigma$. The weight is equivalently the least number of
transpositions whose product can make up $\sigma$. As in equation (A.5) we set
\begin{equation*}
	P_w(n) = \stirlingi{n}{n-w}
\end{equation*}
the number of permutations of weight $w$ in the symmetric  group on $n$ objects.
\bigskip

\noindent
3) 
For disjoint finite sets $A$ and $B$ we define the \emph{cross product} of $\sigma \in \symm(A)$ and $\tau \in \symm(B)$ to be the permutation
\begin{equation*}
	\sigma \times \tau \in \symm(A \cup B)
\end{equation*}
which restricts to $\sigma$ on $A$ and $\tau$ on $B$. The cross product
identifies $\symm(A) \times \symm(B)$ with a subgroup of $\symm(A \cup B)$ which we also denote by $\symm(A) \times \symm(B)$. The cross product extends to products of permutations on three or more pairwise disjoint finite sets. We remark that 
\begin{equation*}
	w(\sigma \times \tau) = w(\sigma) + w(\tau).
\end{equation*}
\bigskip

\noindent
4) For finite sets $A_1,\dots,A_m$ then
\begin{equation*}
	\sum_{w_1+\cdots+w_m = w} P_{w_i}(|A_1|)\cdots P_{w_m}(|A_m|)
\end{equation*}
is the number of $\rho \in \symm(A_1)\times \cdots \times \symm(A_m)$ of weight $w$. Let us write this number as 
\begin{equation*}
	P_w(|A_1|,\dots,|A_m|)
\end{equation*}
for convenience.
\bigskip

\noindent
5) For a positive integer $g$ define $[g] = \{1,\dots,g\}$ and fix some pairwise disjoint finite sets $C_1,\dots,C_g$.  For $I \subseteq [g]$ define
\begin{equation*}
	C_I = \bigcup_{i \in I} C_i
\end{equation*}
and
\begin{equation*}
	C = C_{[g]}= C_1\cup \cdots\cup C_g.
\end{equation*}
 If we have a partition $\mathcal{I}=(I_1,\dots,I_r)$ of $[g]$ into $r$ non-empty sets, then $C_{\mathcal{I}} = (C_{I_1},\dots,C_{I_r})$ is  partition of $C$ into $r$ non-empty sets. This defines 
a subgroup
\begin{align*}
	\symm_{\mathcal{I}}(C) = \symm(C_{I_1})\times \cdots \times \symm(C_{I_r})
\end{align*}
of $\symm(C)$.
\bigskip

\noindent
6) Let $\mathcal{U}_r([g])$ denote the collection of partitions of $[g]$ into $r$ disjoint non-empty parts. Then 
the identity we are trying to prove is the following
\begin{align}
	\sum_{r=1}^g (-1)^{r-1}(r-1)!\sum_{\{I_1,\dots,I_r\} \in \mathcal{U}_r([g])}P_w(|C_{I_1}|,\dots,|C_{I_r}|)=0
\end{align}
for $0 \leq w \leq g-2$. Compare this statement to that of Theorem A.2. The difference of numerical factor
by $r!$ is due to that the subsets in equation (A.11) are ordered. ( The small difference that the ${c_i}$ of (A.11)
need not be integer, but distinct, here they are integer but need not be distinct. )
\bigskip

\noindent
7) To prove (C.1) we consider the left side as a weighted count of permutations of weight $w$ in $\symm(C)$. 
Let $\mathcal{I} = \{I_1,\dots,I_r\}$ be a partition of $[g]$. We say that 
$\mathcal{I}$ \emph{admits} $\sigma \in \symm(C)$ if $\sigma \in \symm_{\mathcal{I}}(C)$, equivalently that $\sigma(C_{I_j}) = C_{I_j}$ for all $j$.
Then
\begin{align*}
	\sum_{r=1}^g (-1)^{r-1}(r-1)!\sum_{\{I_1,\dots,I_r\} \in \mathcal{U}_r([g])}P_w(|C_{I_1}|,\dots,|C_{I_r}|) = \sum_{\sigma \in \symm(C),w(\sigma)=w} \sum_{r=1}^g (-1)^{r-1}(r-1)!\Phi_r(\sigma)
\end{align*}
where $\Phi_r(\sigma)$ is the number of partitions of $[g]$ into $r$ non-empty parts which admit $\sigma$. The  identity ( C.1 ) will follow from the identity
\begin{align}
	\label{eq:x.4}
	\sum_{r=1}^g (-1)^{r-1}(r-1)!\Phi_r(\sigma) = 0
\end{align}
which  I claim holds whenever $w \leq g-2$.
\bigskip

\noindent
8) Given a $\sigma \in \symm(C)$ there is a unique partition of $[g]$,  $\mathcal{J}_\sigma$, such that it
admits $\sigma$, and it is a refinement of any partition that admits $\sigma$.
\bigskip

\noindent
9) Suppose $\mathcal{I}_\sigma$ is a partition into $m$ non-empty subsets. The set of partitions
of $[g]$ with $r$ parts that are refined by $\mathcal{J}_\sigma$ corresponds naturally to the partition of a set of $m$ objects into $r$ parts. Therefore, 
$\Phi_r(\sigma) = \stirlingii{m}{r}$. 
Then eq.(C.2)  is equivalent to
\begin{align}
	\sum_{r=1}^m (-1)^{r-1}(r-1)!\stirlingii{m}{r}=0
\end{align}
\bigskip

\noindent
10) To prove eq.(C.3) we first show $m>1$ if $w \leq g-2$. As was noted in 2), under this condition $\sigma$ can be
realized as the product of $\leq g-2 $ transpositions. The subset in  $\mathcal{J}_\sigma$ that contains $1$ has
$\leq g-1$ elements, since at most $g-2$ elements can be 'connected' to $1$ by the transpositions. So $m>1$.
\bigskip

\noindent
11) We now prove eq.(C.3) is true for $m>1$, and so our result eq.(C.1).
Using the recurrence relation for the $\stirlingii{m}{r}$ gives, for $m\geq 2$

\begin{align*}
	\sum_r (-1)^{r-1}(r-1)!\stirlingii{m}{r} &= \sum_{r} (-1)^{r-1}(r-1)!\left( \stirlingii{m-1}{r-1} + r\stirlingii{m-1}{r} \right) \\
						&= \sum_{r}(-1)^{r-1}(r-1)!\stirlingii{m-1}{r-1}+\sum_{r} (-1)^{r-1}r!\stirlingii{m-1}{r} \\
						&= \sum_{s} (-1)^s s!\stirlingii{m-1}{s} + \sum_{r} (-1)^{r-1}r!\stirlingii{m-1}{r} \\
						&= 0.
\end{align*}

\end{document}